\documentclass[11pt]{article}
\usepackage{amsmath, amsthm, amssymb}
\usepackage{enumerate}
\usepackage{graphicx}
\usepackage{hyperref}
\usepackage[margin=1in]{geometry}
\usepackage{comment}
\numberwithin{equation}{section}

\newtheorem{theorem}{Theorem}[section]
\newtheorem{proposition}[theorem]{Proposition}

\newtheorem*{conjecture*}{Conjecture}
\newtheorem*{OQ*}{Open Question}

\newtheorem{remark}[theorem]{Remark}

\begin{document}
\title{Qualitative analysis of three-wave interaction with periodic boundary condition}

\author{Ze Cheng\footnote{Partially supported by NSF DMS-1405175.}
 ~and Harvey Segur\footnote{Supported by NSF DMS-1107534}
\\
 Department of Applied Mathematics, \\
University of Colorado Boulder 
}
% omit date
\date{\vspace{-5ex}}
% smaller space above title
%\setlength{\droptitle}{-5em}

\maketitle

\begin{abstract}
First, for 3WRI with positive wave energy, we present a regularity theorem for all spatial dimension.
Second, for 3WRI with negative wave energy, we present a class of solution in general spatial dimension that will blow up in finite time. Moreover, a complete classification of spatial uniform solution is given for this particular system.
\end{abstract}

\section{Introduction}
In this paper we study the 3-wave resonance interaction (3WRI) system,
\begin{align}\label{3waveOrigin}
\left\lbrace \begin{array}{cc}
 \partial_{\tau} A_1 +c_1\cdot\nabla A_1 = i \gamma_1\overline{A_2A_3}, \\
 \partial_{\tau} A_2 +c_2\cdot\nabla A_2 = i \gamma_2\overline{A_1A_3}, \\
 \partial_{\tau} A_3 +c_3\cdot\nabla A_3 = i \gamma_3\overline{A_1A_2},   
\end{array}\right. \text{ in }  \Omega, 
\end{align}
with periodic boundary condition, where $\Omega$ is a rectangle domain, 
$$\Omega = \{ x\in \mathbb{R}^n| \ |x_k|<a_k, k=1,\cdots,n \}$$
and all $A_j's$ are are complex amplitude and periodic on $\Omega$.
$\gamma_j=\pm 1$, and $c_j$'s are real non-zero constant vectors.
The derivation of \eqref{3waveOrigin} is rather standard in e.g. nonlinear optics and can be found in e.g. Chap. 4 in \cite{AS81} and \cite{KRB79}. 

Here we are interested in the qualitative property of solution to \eqref{3waveOrigin}. Namely, we are interested in determining if a solution possesses either finite blow up or global existence. Depending on whether $\gamma_j$'s have the same sign, the system \eqref{3waveOrigin} can be classified into two cases, and each case models different physical phenomena.

$\gamma_j$'s \textbf{not having the same sign} corresponds to 3WRI with positive wave energy, which generates from various physical backgrounds e.g. nonlinear optics and water waves. In fact most 3WRI studies in the past focused on this case, and an important  tool is inverse scattering transformation (IST), with which people construct solution and solve the system numerically for some special initial value, e.g. ``separable'' initial value \cite{Kaup80}, or more general initial value \cite{Kaup81}.
%$c_j$'s are assumed to be linearly independent in this case. If $c_j$'s are linearly dependent, \eqref{3waveOrigin} can be reduced to one spatial dimensional 3WRI.  
The 1D-3WRI was approached by Zakharov and Manakov \cite{ZM76}  (and independently  by Kaup \cite{Kaup76}) with inverse scattering transformation, but the boundary condition they considered is that the solution decays sufficiently fast as $|x|\rightarrow\infty$. %Also by implementing IST, Kaup et al. \cite{KRB79} considered the problem with space-independent initial phase. 
A detailed review of the 1-D problem can be found e.g. in \cite{KRB79}. For 3-D problem, Ablowitz and Haberman's work \cite{AH75} leads Cornille \cite{Cornille79} to reformulate the problem into integral equations and 
Kaup gave explicitly general inverse-scattering solution in a series of papers \cite{Kaup80,Kaup81}. See also some recent developments in \cite{GR98}.

However, it seems there has not been any regularity analysis of 3WRI in the case of $\gamma_j$'s \textbf{not} having the same sign. Here, we prove that the system \eqref{3waveOrigin} cannot develop singularity in finite time given that the coefficients $c_j$, $j=1,2,3$, are the same. %If $\gamma_j$'s do \textbf{not} have the same sign, then the constants of motion of system \eqref{3waveOrigin} (see \eqref{constantsMotion}) leads to $L^2$ boundedness of the solution. Thus, by a scaling argument (or ``blow-up'' procedure) we get % can rule out the occurrence of singularity in finite time for $n\leq 2$,
\begin{theorem}\label{regularityPositiveEnergy}
For all space dimension, if initial data is continuously differentiable in space, $\gamma_j$'s do not have the same sign and $c_j$'s $j=1,2,3$ are the same, then the solution to system \eqref{3waveOrigin} exists globally in time.
\end{theorem}

The case that $\gamma_j$'s have \textbf{the same sign} corresponds to 3WRI with negative wave energy, and \eqref{3waveOrigin} becomes \eqref{3waveG}.
%which corresponds to
%\begin{align*}
% k_1+ k_2+ k_3 = 0, \ \ \omega_1+ \omega_2+ \omega_3 = 0.
%\end{align*}
One interesting phenomena about such 3WRI system is that the solution can blow up in finite time. Coppi et al. \cite{CRS69} first found such instability of 3WRI with negative wave energy in plasma physics, which they called ``explosively unstable''.

Heuristically, the transportation wave equation is non-dispersive, and the nonlinear term should enhance the amplitude in a superlinear way. One can compare this to a different system, the 3-D wave equation, where the positive feedback from the nonlinear term needs to race with the dispersive tendency (See a classic paper of F. John \cite{John79}). So, finite-time blow up should not be a surprise in system \eqref{3waveOrigin} when $\gamma_j$'s have the same sign.
However, note that complex amplitude means phase interaction and $c_j$'s can be different, and both cause difficulty in analysis.
This papaer is dedicated to determining if a solution will blow up in finite time given any initial value.

Besides finite-in-time blow-up solution, system \eqref{3waveG} obviously admits global-in-time existing  solution, e.g., two of $A_j$'s are zero and the third is a constant. A natural question is whether \eqref{3waveG} admits other non-trivial globally existing solutions.
In this paper, we show that the system admits globally existing and decaying solution. Moreover, we see that the spatial uniform solution of \eqref{3waveG} can be completely classified in terms of blowup or not. The spatial uniform case of \eqref{3waveG} reduces to an ODE system which has been well studied for decades. Here we state a classification result which may be known to other researchers through various analysis approach. Nevertheless, in Section 5.3.1 we give a analytic proof for convenience of the reader.
\begin{theorem}\label{classificationODE}
A necessary and sufficient condition of finite-in-time blow-up solution to \eqref{3wave} is that the initial condition satisfies one of the following,
\begin{enumerate}
\item Only one of $A_j(0)$'s is zero, i.e., $|A_3(0)|\geq|A_1(0)|>|A_2(0)|=0$; 
%\item  $(\theta_1+\theta_2+\theta_3)(0)=\frac{3\pi}{2}$ and $|A_1(0)|=|A_2(0)|>|A_3(0)|>0$;
\item $(\theta_1+\theta_2+\theta_3)(0)=\frac{3\pi}{2}$, one of $A_j(0)$'s is strictly less than the other two and none is zero, i.e., $|A_3(0)|\geq|A_1(0)|>|A_2(0)|>0$;
\item $(\theta_1+\theta_2+\theta_3)(0)\neq\frac{3\pi}{2}$ (implicitly none of $A_j(0)$'s is zero),
\end{enumerate}
where the indexes $\{1,2,3\}$ allow any permutation and $\theta_j$'s are from $A_j=r_j e^{i\theta_j}$, $j=1,2,3$.
\end{theorem}

\begin{remark}
We think an interesting question is, whether all globally existing solution to \eqref{3waveG} is spatially uniform, and this is so far open.
\end{remark}

We also describe a new class of solutions that blow up in finite time (Theorem \ref{finiteBlowUp}),
\begin{theorem}\label{finiteBlowUp}
Suppose $\theta_j(x,0)=\theta_j(0)$, $j=1,2,3$, and $(\theta_1+\theta_2+\theta_3)(0)=\frac{\pi}{2}$, and  then the solution of \eqref{3waveG} blows up in finite time.
\end{theorem}
%This class of solution can be `inseparable' and hence is new compared to \cite{Kaup80}. 

The mathematical technique used here is elementary, and the results are based on the method of characteristics. We believe that a natural step after this is to implement perturbation theory.

This paper is organized this way. In section 2, we focus on 3WRI with positive wave energy and prove Theorem \ref{regularityPositiveEnergy}. In section 3, we study 3WRI with negative wave energy. In section 3.1, we focus on space-independent case and prove Theorem \ref{classificationODE}. In section 3.2, we consider the general case, and after briefly discussing the well-posedness of \eqref{3waveOrigin} we prove Theorem \ref{finiteBlowUp}.

\section{3WRI with positive wave energy}
Consider system \eqref{3waveOrigin} with $\gamma_1=1$, $\gamma_2,\gamma_3=-1$, and $c_1=c_2=c_3=c$.
%To outline the idea of the proof, we first show a well-known fact that the $L^2$ norm is bounded. Then we assume there exists a singularity at time $t_*$. Then we have a sequence of time approaching $t_*$ at which the maximum norm of the solution approaches infinity. After a rescaling, we get a sequence of bounded solution which converges to a limit solution, and this leads to a contradiction.

Proof of Theorem \ref{regularityPositiveEnergy}. 
Multiply \eqref{3waveOrigin} by $\overline{A_j}$ and take conjugate of the system then multiply by $A_j$, we have
\begin{align*}
\partial_{\tau} |A_1|^2 +c\cdot\nabla |A_1|^2 &= i (A_1A_2A_3+\overline{A_1A_2A_3}), \\
\partial_{\tau} |A_2|^2 +c\cdot\nabla |A_2|^2 &= -i (A_1A_2A_3+\overline{A_1A_2A_3}), \\
\partial_{\tau} |A_3|^2 +c\cdot\nabla |A_3|^2 &= -i (A_1A_2A_3+\overline{A_1A_2A_3}). 
\end{align*}
Hence, we have
\begin{align*}
\partial_{\tau} (|A_1|^2+|A_2|^2) +c\cdot\nabla |A_1|^2 + c\cdot\nabla |A_2|^2 = 0, \\
\partial_{\tau} (|A_1|^2+|A_3|^2) +c\cdot\nabla |A_1|^2 + c\cdot\nabla |A_3|^2 = 0.
\end{align*}
These lead to
\begin{align*}
	|A_1(x+c \tau, \tau)|^2+|A_2(x+c \tau, \tau)|^2 = K_1(x), \\
	|A_1(x+c \tau, \tau)|^2+|A_3(x+c \tau, \tau)|^2 = K_2(x).
\end{align*}
Since the initial data of $A_j$'s are smooth, we know that $K_1(x)$ and $K_2(x)$ must be smooth and bounded. Hence, $|A_j|$, $j=1,2,3$, must be bounded for all time. $\Box$

\begin{remark}
The proof above is obviously valid for general domain and boundary condition. % So, in general we have regularity theorem for 3WRI with positive wave energy and smooth initial data.
\end{remark}

\section{3WRI with negative wave energy}

Now, if $\gamma_j$'s have the same sign, \eqref{3waveOrigin} becomes
%after a simple scaling of $A_j$'s we can move $\gamma_j$'s out of \eqref{3waveOrigin} and get
\begin{align}\label{3waveG}
\left\lbrace \begin{array}{cc}
 \partial_{\tau} A_1 +c_1\cdot\nabla A_1 = i \overline{A_2A_3}, \\
 \partial_{\tau} A_2 +c_2\cdot\nabla A_2 = i \overline{A_1A_3}, \\
 \partial_{\tau} A_3 +c_3\cdot\nabla A_3 = i \overline{A_1A_2},   
\end{array}\right. \text{ in }  \Omega, 
\end{align}
with periodic boundary condition.
The constants of motion are,
\begin{align}\label{constantsG}
K_1 = \int_{\Omega}|A_1|^2-|A_2|^2 dx, \quad K_2 = \int_{\Omega}|A_1|^2-|A_3|^2 dx,
\end{align}
where $K_1, K_2$ are constants. 

\subsection{Spatially uniform case}
There are many classical studies (for reference, see section I. D. of \cite{KRB79}) of space-independent 3WRI which is an ODE system. Here we give a complete classification of space independent solutions of the negative wave energy case of 3WRI,
\begin{align}\label{3wave}
\left\lbrace \begin{array}{cc}
 \partial_{\tau} A_1 = i \overline{A_2A_3}, \\
 \partial_{\tau} A_2 = i \overline{A_1A_3}, \\
 \partial_{\tau} A_3 = i \overline{A_1A_2}.   
\end{array}\right. 
\end{align}
The system \eqref{3wave} admits three constants of the motion:
\begin{align}\label{constants}
\arraycolsep=1.4pt\def\arraystretch{1.3}
\left. \begin{array}{cc}
  K_1 = |A_1|^2-|A_2|^2, \quad K_2 = |A_1|^2-|A_3|^2, \\
  H = A_1A_2A_3+\overline{A_1A_2A_3},
\end{array}\right. 
\end{align} 
where $K_1, K_2$ and $H$ are constants. The first two constants are also called Manley-Rowe relations.

By direct calculation we rewrite \eqref{3wave} as 
\begin{align}
	&\partial_{\tau} r_j^2  = 2  r_1r_2r_3 \sin(\theta_1+\theta_2+\theta_3), \text{ for } j=1,2,3, \label{3waveVar}\\ 
&\left. 
\begin{array}{cc}
	r_1 \partial_{\tau} \theta_1 = r_2r_3\cos(\theta_1+\theta_2+\theta_3), \\
	r_2 \partial_{\tau} \theta_2 = r_1r_3\cos(\theta_1+\theta_2+\theta_3), \\
	r_3 \partial_{\tau} \theta_3 = r_1r_2\cos(\theta_1+\theta_2+\theta_3), 
	\end{array}\right. \label{3waveVarTheta}
\end{align}
and \eqref{constants} becomes
\begin{align}\label{constantsVar}
\arraycolsep=1.4pt\def\arraystretch{1.3}
\left. \begin{array}{cc}
	K_1=r_1^2-r_2^2, \quad  \ K_2=r_1^ 2-r_3^2, \\
	H = 2r_1r_2r_3 \cos(\theta_1+\theta_2+\theta_3).
	\end{array}\right. 
\end{align}

\begin{remark}\label{AiZero}
\eqref{3waveVar}-\eqref{3waveVarTheta} holds only when none of $r_j$'s is zero. This is because $\theta_j$ is not defined when $r_j=0$. 
Actually, we should pay special attention to the situation when $A_j$'s touch zero since \eqref{3wave} does not satisfy Lipchitz condition there. However, we can take advantage of \eqref{constants}. For instance $A_1$ touches zero at $\tau=\tau_0$ and $A_2(\tau_0),A_3(\tau_0)\neq 0$, then $\partial_{\tau}A_j(\tau_0)\neq 0$, i.e., $|A_1|$ will increase and hence $|A_2|,|A_3|$ will increase since $K_1,K_2$ are constant. By \eqref{3waveVar}, we must have $\sin(\theta_1+\theta_2+\theta_3)>0$ and this implies $(\theta_1+\theta_2+\theta_3)(\tau_0+)\in(0,\pi)$. Moreover, we know $H\equiv 0$ and hence $(\theta_1+\theta_2+\theta_3)(\tau)\equiv \frac{\pi}{2}$ for $\tau\in(\tau_0,\infty)$. %Picturing a phase plane of $(r_1,\theta_1+\theta_2+\theta_3)$ in $[\tau_0,\infty)$, the trajectory shoots from origin along the upper vertical axis.
With a little extra effort we can show that the solution must blow up in finite time.

If two of $A_j$'s, say $A_1,A_2$, are zero at $\tau=\tau_0$, the solution of \eqref{3wave} has to be an equilibrium, see the proof of Theorem \ref{globalSolution}.
\end{remark}

\begin{proposition}
If $(\theta_1+\theta_2+\theta_3)(0)\neq\frac{\pi}{2} \text{ or } \frac{3\pi}{2}$, then none of $A_j$'s will touch zero. 
\end{proposition}
Proof. Since $H\neq0$, $A_j\neq 0$. $\Box$

\begin{theorem}\label{globalSolution}
A necessary and sufficient condition for the solution of \eqref{3wave} to exist globally in time is that, the initial data satisfyå either of the following,
\begin{enumerate}
\item[(i)] $|A_3(0)|\geq|A_1(0)|=|A_2(0)|=0$; 
\item[(ii)] $|A_3(0)|\geq|A_1(0)|=|A_2(0)|>0$ and $(\theta_1+\theta_2+\theta_3)(0)=\frac{3\pi}{2}$;
\end{enumerate}
where the indexes $\{1,2,3\}$ allow any permutation.
\end{theorem}

Proof. \textit{We first prove the sufficiency.} 

For (ii) by \eqref{constantsVar} $H,K_1\equiv 0$ and $K_2\leq0$. So, $r_1=r_2$.

Notice that our choice of initial value ensures that all $r_j$, $j=1,2,3$ decay until one of them touches zero. Let's define $\tau_0$ to be the first time that one of $r_j$'s touches zero. If $\tau_0=\infty$ then we have a global decaying solution. So, we assume $\tau_0<\infty$ and to sum up:
\begin{enumerate}
\item $r_j$, $j=1,2,3$ decays in $[0,\tau_0]$;
\item  $(\theta_1+\theta_2+\theta_3)(\tau)=\frac{3\pi}{2}$ in $[0,\tau_0)$;
\end{enumerate}
Let $d_0:=|A_3(0)|-|A_1(0)|\geq0$. \textbf{Case 1:} If $d_0=0$, we can solve the system \eqref{3waveVar} on interval $[0,\tau_0)$ by 
\begin{align}
r_j=\frac{r_j(0)}{1+2r_j(0)\tau}. \label{d0}
\end{align}
Thus $r_j(\tau_0)\neq0$, which contradicts with the assumption on $r_0$.  So, \eqref{d0} can extend to $[0,\infty)$.

\textbf{Case 2:} If $d_0>0$, by \eqref{3waveVar} and the facts $r_1=r_2$, $r_3(\tau)\leq r_3(0)$ on $[0,\tau_0)$ we have
$$\partial_{\tau} r_1 \geq -r_3(0) r_1.$$
Hence by Gr\"onwall's inequality,
\begin{align}\label{d1}
	r_1=r_2\geq r_1(0)e^{-r_3(0)\tau}.
\end{align}
Similar to case 1, $r_j(\tau_0)\neq0$ contradicts the assumption on $r_0$, so \eqref{d1} can extend to $[0,\infty)$.

For (i), the initial data obviously admits an equilibrium solution, and we claim that this is the only solution (this is not obvious since we lose Lipchitz condition of \eqref{3wave} if some $A_j$'s are zero). Suppose a non-constant solution $A_1(\tau),A_2(\tau),A_3(\tau)$ with initial condition $|A_3(0)|\geq|A_1(0)|=|A_2(0)|=0$ (so $H\equiv 0$) and $|A_3|\geq|A_1|=|A_2|>0$ at $\tau=\tau_0$ (since $K_1\equiv 0$, $r_1\equiv r_2$). Then if we change the time direction of \eqref{3wave} by letting $\tau=\tau_0-\eta$ and $B_j(\eta)=A_j(\tau_0-\eta)$ (and corresponding new $r_j$ and $\theta_j$), $j=1,2,3$, then we get a new system,
\begin{align}
\left\lbrace \begin{array}{cc}
 \partial_{\eta} B_1 = -i \overline{B_2B_3}, \\
 \partial_{\eta} B_2 = -i \overline{B_1B_3}, \\
 \partial_{\eta} B_3 = -i \overline{B_1B_2},   
\end{array}\right. 
\end{align}
with initial value $|B_3(0)|\geq |B_2(0)|=|B_1(0)|> 0$ at $\eta = 0$. From $H\equiv 0$ ($H,K_1,K_2$ do not change) and $r_j,j=1,2,3$, decaying for some interval, we see that $\theta_1+\theta_2+\theta_3=\frac{\pi}{2}$ for that interval. So, similar to the proof for (ii), it has a unique solution that decays to a limit $(0,0,A_3(\tau=0))$ as $\eta\rightarrow \infty$ and never touches zero in finite time. This contradicts with $(B_1(\tau_0),B_2(\tau_0),B_3(\tau_0))=(A_1(0),A_2(0),A_3(0))=(0,0,A_3(0))$.

\quad

\textit{Now we turn to the necessity.}  If the initial data does not satisfy either condition in Theorem \ref{globalSolution}, we have the following cases,
\begin{enumerate}
\item Only one of $A_j(0)$'s is zero, i.e., $|A_3(0)|\geq|A_1(0)|>|A_2(0)|=0$; 
%\item  $(\theta_1+\theta_2+\theta_3)(0)=\frac{3\pi}{2}$ and $|A_1(0)|=|A_2(0)|>|A_3(0)|>0$;
\item $(\theta_1+\theta_2+\theta_3)(0)=\frac{3\pi}{2}$, one of $A_j(0)$'s is strictly less than the other two and none is zero, i.e., $|A_3(0)|\geq|A_1(0)|>|A_2(0)|>0$;
\item $(\theta_1+\theta_2+\theta_3)(0)\neq\frac{3\pi}{2}$ (implicitly none of $A_j(0)$'s is zero).
\end{enumerate}
Again, the above indexes $\{1,2,3\}$ allow any permutation. We are going to show that solutions to the above initial conditions will blow up in finite time.

For \textbf{case 1}, by Remark \ref{AiZero}, we know $(\theta_1+\theta_2+\theta_3)(\tau)\equiv \frac{\pi}{2}$ for $\tau\in(0,\infty)$, and by \eqref{3waveVar} and \eqref{constantsVar} we have,
$$\partial_{\tau} r_1^2  = 2  r_1\sqrt{r_1^2-K_1}\sqrt{r_1^2-K_2},$$
where $K_1,K_2<0$. Hence it is easy to see $r_1$ blows up in finite time.

For \textbf{case 2}, by \eqref{3waveVar} we know that $A_2$ will touch zero first in finite time, call this time as $\tau_0$. Then $\partial_{\tau}A_2(\tau_0)\neq 0$ and by Remark \ref{AiZero} we know $(\theta_1+\theta_2+\theta_3)(\tau)=\frac{\pi}{2}$ for $\tau>\tau_0$. This implies a finite-in-time blow up.

For \textbf{case 3}, let $\Theta=\theta_1+\theta_2+\theta_3$. If $\Theta(0)\in(-\frac{\pi}{2},\frac{\pi}{2}]$ (we change $\frac{3\pi}{2}$ to $-\frac{\pi}{2}$ for convenience), 
\textbf{step 1.}  we show that
$\Theta$ increases and approaches $\frac{\pi}{2}$.
%By \ref{constantsVar} this means $r_1r_2r_3$ going to infinity.
By \eqref{3waveVarTheta} we have
\begin{align*}
	 \partial_{\tau} \Theta = (\frac{r_2r_3}{r_1}+\frac{r_1r_3}{r_2}+\frac{r_1r_2}{r_3})\cos\Theta.
\end{align*}
By \ref{constantsVar} we have $|r_1r_2r_3|\geq \frac{|H|}{2}$. Suppose $r_j>\eta_1$, $j=1,2,3$ for some $\eta_1>0$, then there exists some constant $\eta_0>0$ such that $(\frac{r_2r_3}{r_1}+\frac{r_1r_3}{r_2}+\frac{r_1r_2}{r_3})>\eta_0$. Now without loss of generality suppose $r_1\rightarrow 0$, then $r_2r_3\rightarrow \infty$, hence $\frac{r_2r_3}{r_1}\rightarrow\infty$. In both cases, we see that $\cos\Theta$ is positive since $\Theta\in(-\frac{\pi}{2},\frac{\pi}{2}]$, and hence $\Theta$ will approach $\frac{\pi}{2}$. 
%Otherwise $\cos\Theta>\eta_2$ for some constant $\eta_2>0$, then $\partial_{\tau} \Theta>\eta_0\eta_2$ and in finite time $\Theta$ will get larger than $\cos^{-1}\eta_2$ ($\Theta$ must get larger than 0), a contradiction.

\textbf{Step 2.} There exists $\tau_1>0$ such that $\Theta>\frac{\pi}{4}$ for $\tau>\tau_1$. Then by \eqref{3waveVar} for  $\tau>\tau_1$
\begin{align*}
\partial_{\tau} r_1  &\geq Cr_2r_3 \\
				&\geq C\sqrt{r_1^2-K_1}\sqrt{r_1^2-K_2},
\end{align*}
where the second inequality is due to \eqref{constantsVar}. So a finite blow up in $r_1$ can be easily derived.

Similarly we have a finite blow up for $\Theta(0)\in(\frac{\pi}{2},\frac{3\pi}{2})$.
$\Box$

Hence Theorem \ref{classificationODE} directly follows from Theorem \ref{globalSolution}.

\subsection{General case}
First we briefly discuss the well-posedness of \eqref{3waveOrigin}. The local existence and uniqueness is guaranteed by method of characteristics. Also, the spatially periodic boundary condition is directly satisfied if initial data are periodic, since all $c_j$'s are constant vectors.

Now we consider \eqref{3waveG},
and by direct calculation we rewrite the system as 
\begin{align}
	&\partial_{\tau} r_j^2  +c_j\cdot\nabla r_j^2 = 2  r_1r_2r_3 \sin(\theta_1+\theta_2+\theta_3), \text{ for } j=1,2,3, \label{3waveGVar}\\ 
&\left. \begin{array}{cc}
	r_1 (\partial_{\tau} \theta_1 +c_1\cdot\nabla\theta_1) = r_2r_3\cos(\theta_1+\theta_2+\theta_3), \\
	r_2 (\partial_{\tau} \theta_2 +c_2\cdot\nabla\theta_2) = r_1r_3\cos(\theta_1+\theta_2+\theta_3), \\
	r_3 (\partial_{\tau} \theta_3 +c_3\cdot\nabla\theta_3) = r_1r_2\cos(\theta_1+\theta_2+\theta_3), 
	\end{array}\right. \label{3waveGVarTheta}
\end{align}
%and \eqref{constantsG} becomes
%\begin{align}\label{constantsGVar}
% const=\displaystyle\int_{\Omega}r_1^2-r_2^2 dx, \quad  \ const=\displaystyle\int_{\Omega}r_1^ 2-r_3^2 dx,
%\end{align}
To outline the proof, if initially $A_j$ is lined up in the same direction at each point, and the summation of $\theta_j$'s equal $\frac{\pi}{2}$, then all $\theta_j$'s preserve at all time. Hence \eqref{3waveGVarTheta} is moved out of the equations. Then we are left to analyze \eqref{3waveGVar}.

Proof of Theorem \ref{finiteBlowUp}. \textbf{Step 1.} We show that $(\theta_1+\theta_2+\theta_3)(x,\tau)\equiv \frac{\pi}{2}$.

Note that implicitly, the initial condition guarantees that $A_j(x,0)$'s are nowhere zero and bounded. So, the short time existence and uniqueness of solution to \eqref{3waveGVar} is provided by method of characteristics. Hence, $\theta_j(x,\tau)\equiv\theta_j(0)$ and $(\theta_1+\theta_2+\theta_3)(x,\tau)\equiv \frac{\pi}{2}$.

\textbf{Step 2.} \eqref{3waveGVar} becomes
\begin{align*}
\left\lbrace \begin{array}{cc}
	\partial_{\tau} r_1  +c_1\cdot\nabla r_1 =  r_2r_3, \\
	\partial_{\tau} r_2  +c_2\cdot\nabla r_2 =  r_1r_3, \\
	\partial_{\tau} r_3  +c_3\cdot\nabla r_3 =  r_1r_2,
	\end{array}\right. \text{ in } \Omega.
\end{align*}
By method of characteristics, the corresponding characteristic equation for $r_1$ is 
\begin{align}\label{characteristicsR}
\arraycolsep=1.4pt\def\arraystretch{1.3}
\left\lbrace \begin{array}{ll}
	\frac{dr_1}{d\tau}  &=  r_2r_3, \\
	\frac{dx}{d\tau}  &= c_1,
	\end{array}\right.
\end{align}
with initial data $x(0)=\xi$ and $r_1(\xi,0)=\phi_1(\xi)$. 
%The solution is
%\begin{align*}
%\arraycolsep=1.4pt\def\arraystretch{1.3}
%\left\lbrace \begin{array}{ll}
%	 r_1(\xi, \tau)&=  \phi_1(\xi) +\int_0^{\tau} r_2(\xi, s)r_3(\xi, s) ds, \\
%	 x(\xi,\tau) &= c_1 \tau + \xi,
%	\end{array}\right.
%\end{align*}
%or written as
%\begin{align*}
%	r_1(x, \tau)&=  \phi_1(x-\tau c_1)+\int_0^{\tau} r_2(x-\tau c_1, s)r_3(x-\tau c_1, s) ds.
%\end{align*}
Since $r_2,r_3\geq 0$, for fixed $\xi$, $r_1(\xi, \tau)$ is monotone increasing as $\tau$ increases. We say $r_1$ is ``monotone increasing'' along characteristics. Similarly, $r_2,r_3$ are monotone increasing along their characteristics $x=c_2\tau+\xi$ and $x=c_3\tau+\xi$.

Since $r_j(0)$'s are strictly positive in $\Omega$, the characteristic equations also imply that $r_j(x,\tau)>0$ on $\overline{\Omega}$ for all $\tau>0$. This eliminates the possible zero at the boundary of $\Omega$. Hence after any $\tau=\tau_0>0$, $r_j$'s must have a positive infimum in $\overline{\Omega}$. Suppose
\begin{align*}
	f(\tau) = \min_{x\in\Omega} \left\lbrace r_1(x,\tau),r_2(x,\tau),r_3(x,\tau)\right\rbrace ,
\end{align*}
then $f(\tau)>0$ for $\tau>\tau_0>0$. We \textbf{claim} that $f(\tau)$ is Lipschitz, and hence from \eqref{characteristicsR} we see
\begin{align*}
\frac{df}{d\tau} \geq f^2,
\end{align*}
hold in weak sense. This implies that $f$ blows up in finite time, since for any positive test function $v\in C_0^1(\mathbb{R})$ and if $g$ is the solution of
\begin{align*}
\arraycolsep=1.4pt\def\arraystretch{1.3}
\left\lbrace \begin{array}{ll}
	\displaystyle\frac{dg}{d\tau} = g^2, \\
     g(\tau_0)=C_0, \text{ and } 0<C_0<f(\tau_0),
	\end{array}\right.
\end{align*}
then $w=f-g$ satisfies
\begin{align*}
  -\int_{\tau_0}^{T} w v' ds \geq \int_{\tau_0}^{T} \beta(s) w(s) v(s)ds, 
\end{align*}
where $\beta(\tau)$ is non-negative (since $f,g$ are positive) and suppose $[\tau_0,T]$ is the support of $v$. By choosing suitable $v$ we see that
\begin{align*}
		w(T) \geq C\int_{\tau_0}^{T} \beta(s) w(s) ds,
\end{align*}
where $v(s)$ on the right is absolved in $\beta(s)$. Hence by Gr\"onwall's inequality $w\geq 0$. Then $g$ blowing up in finite time implies that $f$ blows up in finite time. 

We finish the proof by proving the claim. Let $\bar{r}_1(\tau)=\min_{x\in\Omega}r_1(x,\tau) = r_1(x_*,\tau)$,. Note that $r_1$ is monotone increasing along characteristics, then $\bar{r}_1(\tau)$ must be monotone increasing. Hence,
\begin{align*}
0\leq \frac{\bar{r}_1(\tau+\Delta \tau)-\bar{r}_1(\tau)}{\Delta \tau} \leq \frac{r_1(x_*,\tau+\Delta \tau)-r_1(x_*,\tau)}{\Delta \tau} \leq |r_1|_{C^1}.
\end{align*}
So $\bar{r}_1(\tau)$ is Lipschitz. Similarly $\bar{r}_2,\bar{r}_3$ are Lipschitz. Finally, $f=\min\{\bar{r}_1,\bar{r}_2,\bar{r}_3\}$ is Lipschitz.
$\Box$

\begin{remark}
Although $\theta_j$'s are constant, $r_j$'s vary in space. So, for the future perturbation work, one only needs to perturb $\theta_j$'s.
\end{remark}

\section*{Acknowledgment}
%ACknowlegement I thank Professor Harvey Segur for introducing to me the problem of 3-wave interaction and all the helpful comments and critiques on this work.
The authors thank Professor Congming Li for his comments and insight of the problem. Also, thank Xinshuo Yang for helpful discussion on this paper.

% To make the spacing of the reference smaller
\providecommand{\bysame}{\leavevmode\hbox to3em{\hrulefill}\thinspace}
\providecommand{\MR}{\relax\ifhmode\unskip\space\fi MR }
\providecommand{\MRhref}[2]{%
  \href{http://www.ams.org/mathscinet-getitem?mr=#1}{#2}
}
\providecommand{\href}[2]{#2}

%\bibliographystyle{plain}
%\bibliography{mainbib}

\noindent
ze.cheng@colorado.edu \\
segur@colorado.edu

\end{document}